\documentclass[12pt]{article}
\usepackage{fullpage}

\usepackage{amsmath,amssymb,amsthm,amsfonts}
\usepackage{enumerate}
\usepackage{picinpar}
\usepackage{amsfonts}
\usepackage{layout}
\usepackage{amsthm, verbatim}
\usepackage{graphicx}
\theoremstyle{plain}

\newtheorem{theorem}{Theorem}[section]
\newtheorem*{theorem*}{Theorem}

\newtheorem{lemma}[theorem]{Lemma}

\newtheorem*{corollary*}{Corollary}

\newtheorem*{conjecture*}{Conjecture}

\theoremstyle{remark}

\theoremstyle{definition}

\newtheorem{example}[theorem]{Example}
\numberwithin{figure}{section}

\newcounter{fig}

\newcommand{\lskip}{\phantom{.}}

\def\R{\mathbb R}

\def\P{\mathbb P}

\def\xl{x^{*}}
\def\E{\mathbb E}

\def\ch2{5,\hspace{-.05in}10\hspace{-.05in}-\hspace{-.05in}CH_2\hspace{-.05in}
  -\hspace{-.05in}THF}

\def\sh2{5,\hspace{-.02in}10\hspace{-.02in}-\hspace{-.02in}CH_2\hspace{-.02in}
  -\hspace{-.02in}THF}

\def\h={5,\hspace{-.05in}10\hspace{-.05in}-\hspace{-.05in}CH
  \hspace{-.05in} =\hspace{-.05in}THF}

\def\s={5,\hspace{-.02in}10\hspace{-.02in}-\hspace{-.02in}CH
  \hspace{-.02in} =\hspace{-.02in}THF}

\def\f10{10f\hspace{-.05in}-\hspace{-.05in}THF}

\def\sf10{10f\hspace{-.02in}-\hspace{-.02in}THF}

 \def\Rom#1{\uppercase\expandafter{\romannumeral
    #1}}

\usepackage{times}

\begin{document}

\pagestyle{plain}

\begin{center} 
{\Large \bf Propagation of Fluctuations in Biochemical Systems, I: Linear SSC Networks }\\[2in]

David Anderson,$^{1}$, Jonathan Mattingly,$^{1}$, H. Frederik Nijhout,$^{2}$  Michael Reed,$^{1*}$

\end{center}

\vspace{3in}

\vspace{1in} 
\hrule 

\vspace{.125in}
\noindent {\small $^{1}$ Department of Mathematics, Duke University,
  Durham, NC 27708; $^{2}$  Department of Biology, Duke University,
  Durham, NC 27708;  

\vspace{.125in}
\noindent $^*$ corresponding author, email:reed@math.duke.edu,
phone:919-660-2808, Fax:919-660-2821}

\newpage
\begin{center}
{\Large \bf Abstract}
\end{center}

We investigate the propagation of random fluctuations through
biochemical networks in which the concentrations of species are large
enough so that the  unperturbed problem is well-described by ordinary
differential equation. We characterize the   behavior of variance as
fluctuations propagate down chains, study the  
effect of side chains and feedback loops, and investigate the
asymptotic behavior as one rate constant gets large. We also describe
how the ideas can be applied to the study of methionine metabolism.

\newpage

 \section{Introduction.}

 There are two different natural contexts in which stochastic dynamics
 arises in the study of biochemical reaction networks. In the first,
 the stochastic chemical dynamics arises from the randomness inherent
 in the formation and breaking of chemical bonds. This ``intrinsic
 stochasticity'' is particularly relevant when the numbers of
 molecules are small such as in gene transcription and small gene
 regulatory networks where the mean concentrations no longer
 faithfully model the chemical dynamics. There is a large literature
 in this field beginning with \cite{delbruck}, including
 \cite{kurtz72}, \cite{gill}, and recently exemplified by
 \cite{othmer}\cite{kurtz05}. In this setting, one typically assumes
 that the reaction system is described by a Poisson process that
 models individual discrete chemical reactions. One then derives a
 partial differential equation for the time evolution of concentration
 densities. As all species have their own intrinsic stochasticity,
 this partial differential equation is parabolic with a uniformly
 elliptic generator.

 In the second context, which is our focus here, one wants to
 investigate the response of a large biochemical system to external
 excitation. It is natural and theoretically useful to consider
 stochastic excitations and to study the emergent properties of the
 network as the random fluctuations propagate through the system. Here
 the randomness is a tool used to study the out-of-equilibrium
 dynamics of the biochemical system. In this setting, we assume that
 the concentrations are large enough so that the unperturbed dynamics
 is faithfully modeled by ordinary differential equations. Typically,
 one in interested in perturbing a single (or small number of)
 input(s) with white noise. Hence, the perturbed problem becomes a
 stochastic differential equation with a hypoelliptic generator.

The central biological goal driving our work is to understand the
behavior of biochemical systems in cells, which {\em in vivo} are
exceptionally large and complicated. A metabolite can be the substrate
for many different enzymes and participate in apparently unrelated
reactions.  Individual reactions usually have nonlinear kinetics
catalyzed by enzymes that are themselves inhibited or excited by
products or distant substrates in the network.  Cells and tissues
differ because the genes that code for certain enzymes have tissue
specific expression patterns and biochemical substrates themselves
also influence gene expression.  Further, each cell's environment, its
inputs and outputs, and its internal state (e.g. stage of cell cycle)
are not constant but vary in time. This continual variation affects
both the concentrations of substrates and the expression of genes that
catalyze particular reactions.  Thus, the gene-biochemical network
should not be viewed as a fixed object but as one that is continuously
changing.

For each signal, either external or internal, that causes a particular
cell to dramatically change its operation, there are two natural
questions. First, how does the gene-biochemical network respond to
accomplish the change? Second, how does the network enable the cell to
maintain homeostasis in all its other operations despite the change?
One would like to understand the structural and kinetic principles
that allow the network to accomplish both tasks simultaneously. We
take two distinct approaches to this biological goal. First, we study
how fluctuations propagate through relatively simple systems. We are
interested in discovering how different network geometries magnify or
suppress fluctuations since this may give clues to why biochemical
networks look the way they do.  Secondly, we apply fluctuations to
{\em in silico} representations of specific biological networks. By
observing how fluctuations propagate we can identify reactions or
subsystems that are buffered against such fluctuations, i.e. are
homeostatic. Then, through {\em in silico} experimentation (e.g.
removing particular reactions), we can take the system apart piece by
piece to discover the regulatory mechanisms that give rise to the
homeostasis.

In this paper, we develop fluctuation theory for chemical reaction
systems for which each complex (in Feinberg's terminology,
\cite{feinberg1}) consists of a single chemical species and the
kinetics are mass action. Thus the corresponding differential
equations are linear, so we refer to such networks as linear SSC
(single species complexes) networks.  Because of the linearity, the
technical difficulties involved in studying the associated stochastic
processes are minimized.  Thus, linear SSC systems are an excellent
arena for investigating the effects of network geometry on the
propagation, magnification, and suppression of fluctuations. The
principles discovered then become the natural goal for generalization
to nonlinear settings \cite{anderson2}.

To see the kinds of questions we want to ask, consider a simple chain
with a side branch.  The chemical species are $X_1, \ldots, X_n, X_s$;
the corresponding concentrations are denoted by $x_1, \ldots, x_n,
x_s$.

{\centering
\begin{picture}(430,90)

\put(75,60){\vector(1,0){50}}
\put(140,60){\makebox(0,0){$X_1$}}
\put(197,53){\vector(0,-1){30}}
\put(191,23){\vector(0,1){30}}
\put(152,60){\vector(1,0){30}}
\put(100,68){\makebox(0,0){$I + \sigma dB(t)$}}
\put(177,38){\makebox(0,0){$k_{s,2}$}}
\put(210,38){\makebox(0,0){$L$}}
\put(168,68){\makebox(0,0){$k_1$}}
\put(194,60){\makebox(0,0){$X_2$}}
\put(194,10){\makebox(0,0){$X_s$}}
\put(206,60){\vector(1,0){30}}
\put(220,68){\makebox(0,0){$k_2$}}
\put(272,68){\makebox(0,0){$k_3$}}
\put(247,60){\makebox(0,0){$X_3$}}
\put(260,60){\vector(1,0){30}}
\put(301,60){\makebox(0,0){$\cdots$}}
\put(320,60){\makebox(0,0){$X_n$}}
\put(332,60){\vector(1,0){40}}
\put(350,68){\makebox(0,0){$k_n$}}
\put(390,30){\makebox(0,0){$.$}}
\end{picture}}

\vspace{.125in}
\noindent The chain has a constant input $I$, which is perturbed by
some random process, in this case, white noise.  If the input is
fluctuating, then each of the concentrations will fluctuate as will
the fluxes, $k_ix_i$. Suppose the side chain is absent. Then, will the
variations of the fluxes increase, decrease, or stay the same as we
move down the chain? Does the answer depend on the rate constants
$k_i$?  If the side chain is present, does it affect the variances of
the fluxes on the chain? If so, what is the effect of the size of $L$.

   The chemical reaction diagram corresponds to set of
differential equations for the concentrations and, similarly, the
diagram with stochastic forcing corresponds to a system of stochastic
differential equations (SDEs):

\vspace{-.25in}
\begin{align*}
  d x_1 &= (I - k_1x_1) dt + \sigma  dB(t)\\
  \dot x_2 &= k_1x_1 - Lx_2 - k_2 x_2 + k_{s,2}x_s\\
  \dot x_3 &= k_2x_2 - k_3x_3\\
  &\ \ \vdots
\end{align*}
These SDEs in turn give rise to a stochastic process on the state
space $\R^{n+1}.$ We prove that this stochastic process has a unique
stationary measure. Intuitively, this means that at large times the
joint distribution of values of the concentrations becomes independent
of the initial condition and independent of time. That is, the
statistics converge to an equilibrium distribution. The variances of
the concentrations referred to above are the variances of the marginal
distributions of this measure. We prove the existence of the
stationary measure for linear SSC systems in Section 2.2. In Section 3
we study the propagation of fluctuations in chains.  In Section 4 we
study the effects of side reaction systems, and feedback loops.  In
Section 5 we ask what happens to variances in the asymptotic limit as
one of the rate constants goes to $\infty,$ corresponding to a very
fast reaction. In Section 6 we show how to use the fluctuation theory
ideas to investigate methionine metabolism.

It is important to note that our, goals, methods and results are
different from those in classical biochemical control theory
\cite{kacser},\cite{crabtree},\cite{heinrich},\cite{westerhoff}.  In
that theory one takes a system at a fixed steady state, makes a small
perturbation in a parameter (perhaps an input), and allows the system
to relax to a new steady state. By comparing the new value of a
variable (a concentration or flux) to the old value, one computes the
percentage change of the variable per unit percentage change in the
parameter. Technically, one is computing a partial derivative. This
kind of sensitivity analysis gives good information about local,
linearized behavior near the initial steady state.  By contrast, we
are concerned with responses to large scale fluctuations in inputs.
Technically, this means computing properties of the distribution of
each concentration or flux from properties of the stationary measure.

It is true that the classical biochemical control theory can be made
``stochastic'' in the following way. Suppose that the system has input
$I$ and is at steady-state.  Consider the same system with input $I
+\eta$, where $\eta$ is a random variable drawn from some density. For
each $\eta$ we let the system relax to steady state and measure the
value, $v$, of some concentration or flux. $v$ is a random variable
and comparing it's variance to the variance of $\eta$ gives
information about how much {\bf the steady state value of} $v$ changes
as $\eta$ changes.  However, this modified biochemical control theory
often gives completely different answers from the fluctuation theory
that we are developing and the differences are biologically
significant.  Consider the chain (without the side chain) in the
example above. If the input is $I +\eta$, then, at steady state, the
flux $k_nx_n$ must equal $I +\eta$, so $Var(k_nx_n) = Var(\eta);$ thus
the variance remains constant down the chain.  By contrast, we will
see below that in our fluctuation theory, under a variety of
reasonable assumptions, that the variances of the fluxes {\em
  decrease} as one proceeds down the chain.  This result is
interesting from a biological point of view because it says that one
way to stabilize the flux out of a chain (i.e. small variance) is to
have many intervening biochemical steps between the input and the
output.

\section{SSC networks and the stationary measure}

In this section we introduce the class of chemical reaction systems
that we will study and prove the existence of a stationary measure.

\subsection{SSC systems with mass actions kinetics}

Throughout we use the terminology introduced by Horn, Jackson, and
Feinberg \cite{horn}\cite{feinberg1}. Let $m$ be the number of
chemical species. We shall study chemical reaction systems such that
each complex contains a single chemical species and refer to such
systems as {\em SSC networks}. In the sequel, we use only the statements in Lemma 2.3.

\begin{lemma}[Deficiency of SSC networks]
An SSC network  has  deficiency  zero.
\label{linearweakreverse}
\end{lemma}
\begin{proof}
  Suppose the network has a single linkage class and let $S$ denote
  the stoichiometric subspace.  Choose any reaction in the network,
  $X_i \rightarrow X_j$. Here we have two complexes and one reaction
  vector in $S$.  Thus, if there are no other complexes, we are done.
  Because the diagram has one linkage class, if there are other
  complexes, then there must be one, call it $X_k$, with an arrow to
  or from either $X_i$ or $X_j$. This adds one complex and one
  dimension to $S$ since $X_k$ is not a linear combination of $X_i$
  and $X_j$. Continuing in this manner until we have exhausted all the
  complexes, we see that the number of complexes is one greater than
  $dim\{S\}.$ Since there is one linkage class the deficiency of the
  network is zero. The case where there is more than one linkage class
  follows easily because the reaction vectors in different linkage
  classes are orthogonal.
\end{proof}

\lskip

We will concentrate on SSC networks containing the zero complex that
have one linkage class.

\begin{lemma}[Dimension of $S$]
  In an SSC system containing the zero complex with one linkage class,
  $dim\{S\} = m$.
\label{defcorollary}
\end{lemma}  
\begin{proof}
  Since the network contains the zero complex, the number of
  complexes, $n$, is one greater than the number of species, $m$. If
  $s = dim\{S\}$, then, by Lemma \ref{linearweakreverse}, $0 = n - s -
  1 = m - s$, so $s=m.$
\end{proof}

\lskip

We assume mass action kinetics so the differential equations governing
the system are linear:
\begin{equation}
  \dot x(t) = A x(t) + I,
  \label{lineareq}
\end{equation}
where $A \in {\mathbb R}^{m \times m}$ and $x(t),I \in {\mathbb R}^m$.
The matrix $A$ is the matrix of rate constants for the system and the
vector $I$ represents any constant flow into the species of the system
from the zero complex. Thus the components of $I$ are non-negative. We
denote the open positive orthant and its closure by ${\mathbb
  R}^m_{>0}$ and $\R^m_{\geq 0}$, respectively.
\begin{lemma}
  If a linear SSC system is weakly reversible and contains the zero
  complex, then
  \begin{itemize}
  \item[(a)] The differential equations \eqref{lineareq} have a unique
    equilibrium which is globally asymptotically stable and contained
    in ${\mathbb R}^m_{>0}$.
  \item [(b)] The eigenvalues of the matrix of rate constants, $A$,
    have strictly negative real parts.
  \item [(c)] For all vectors $v \in \R^m_{\geq 0}$, we have
    $e^{At}v\cdot e_j \geq 0$.
  \end{itemize}
  \label{linearstable}
\end{lemma}
\begin{proof}
  Part (a) is a special case of the zero deficiency theorem
  \cite{feinberg1}. Since $A$ is the Jacobian at the equilibrium
  point, (b) follows from (a) and linearity. (c) holds because
  $\R^m_{\geq 0}$ is invariant under the flow of the differential
  equation.
\end{proof}

\subsection{The Stationary Measure.}

Consider the following weakly reversible SSC system with mass action
kinetics, input vector $I$, and matrix of rate constants $A$ perturbed
by a mean zero, finite variance stationary stochastic process $\xi(t)$:
\begin{equation}
  \left\{
    \begin{aligned}
      \dot x(t) &=  Ax(t) + I \; + \;  \xi(t)\;,\\
      x(0)&=x_0\;.
    \end{aligned}
  \right.
  \label{linearsde}
\end{equation}
From this definition and the stationarity of $\xi(t)$ one easily sees
that \eqref{linearsde} generates a time-homogeneous Markov process.
\begin{theorem}  
  The process $\xl(t)=\xl(t,\xi)$ defined by
  \begin{equation}
    \label{eq:xStationary}
    \xl(t,\xi) = \int_{-\infty}^t e^{A(t-s)}I \ ds + \int_{-\infty}^t e^{A(t-s)} \xi_s \ ds
  \end{equation}
  is a stationary solution to \eqref{linearsde}. Furthermore given any
  initial condition $x_0$, if $x(t,x_0,\xi)$ is a solution to equation
  \eqref{linearsde} then $x(t,x_0,\xi)$ converges to $\xl(t,\xi)$ as
  $t \rightarrow \infty$ in that
  \begin{equation*}
    \E |x(t,x_0,\xi) - \xl(t,\xi)|^2 \rightarrow 0 \quad\text{ as}\quad t \rightarrow \infty\;. 
  \end{equation*}
  \label{convergence}
\end{theorem}
\begin{proof}
  Observe that for any $t,\tau \in \R$,
  \begin{align*}
    \xl(t+\tau) &= \int_{-\infty}^{t+\tau} e^{A(t + \tau -s)}I \ ds +
    \int_{-\infty}^{t+\tau} e^{A(t+\tau - s)}\xi(s) \ ds \notag\\
    &=  \int_{-\infty}^t e^{A(t-s)}I  \ ds + \int_{-\infty}^t e^{A(t-s)}
    \xi(s+\tau) \ ds\;. 
  \end{align*}
  This can be written succinctly as
  \begin{equation}\label{eq:fiberShift}
    (\theta_\tau \xl)(t,\xi)= \xl(t, \theta_\tau\xi)
  \end{equation}
  where the shift $\theta_t$  is defined by $(\theta_t f)(s) = f(t+s)$ for
  all $s,t \in \mathbb R$ and functions $f$ on $\R$. 
  Hence for any $t_1 \leq \cdots \leq t_n$, 
  \begin{equation*}
    \big(\xl(\tau+t_1,\xi),\cdots, \xl(\tau+t_n,\xi)\big)=
    \big(\xl(t_1,\theta_\tau\xi),\cdots, \xl(t_n,\theta_\tau\xi)\big) \;.
  \end{equation*}
  Since $\xi$ is a stationary process, the distribution of the right
  hand side is independent of $\tau$ which proves that $\xl$ is
  stationary. Clearly, $\xl(t,\xi)$ is a solution in that
  $x(t,\xl(0,\xi),\xi)=\xl(t,\xi)$.

  We now turn to convergence. It follows from Lemma 2.3(b) that there
  are constants $\alpha, M >0$ such that $\|e^{At}\| < M e^{-\alpha t}$
  for all $t>0$.  Subtracting the solution of (2),
  \begin{equation} 
    x(t,x_0,\xi) = e^{At}x_0 + \int_{0}^t e^{A(t-s)}I  \ ds + \int_{0}^t e^{A(t-s)} \xi_s \ ds,
  \end{equation}
  from $\xl(t)$, squaring, and taking expected values gives, 
  \begin{align*}
    \E |x(t,x_0,\xi) - \xl(t,\xi)|^2 \leq& 3\|e^{At}\|^2 |x_0|^2 + 3\E
    \left | \int_{-\infty}^0 e^{A(t-s)}I ds \right|^2 + 3\E \left |
      \int_{-\infty}^0 e^{A(t-s)} \xi_s  ds \right|^2\\
    \leq& 3M^2 |x_0|^2 e^{-2\alpha t} + \frac{3M^2
      |I|^2}{\alpha^2}e^{-2
      \alpha t}\\
    & \hspace{.65in} + 3 \E \left( \int_{-\infty}^0 \|e^{A(t-s)}\| ds
    \right)\left( \int_{-\infty}^0 \|e^{A(t-s)}\| |\xi_s|^2 ds \right)\\
    \leq& 3M^2 |x_0|^2 e^{-2\alpha t} + \frac{3M^2
      |I|^2}{\alpha^2}e^{-2 \alpha t} + \frac{3M^2}{\alpha^2}
    e^{-2\alpha t}Var(\xi).
  \end{align*}
  Thus, $\E |x(t,x_0,\xi) - \xl(t,\xi)|^2 \rightarrow 0 \quad\text{ as}\quad t \rightarrow \infty\;.$ 
\end{proof}

\noindent {\bf Remark.} If one takes expectations on both sides of
equation (5), one sees immediately that the model is consistent in the
mean, that is, the mean of the perturbed problem is equal to the
solution of the unperturbed problem.

\vspace{.125in}

If instead of random perturbations given by the vector $\xi_t$ we had
allowed the system to be perturbed by independent white noise
processes, we arrive at the following system of It\^o stochastic
differential equations:
\begin{equation}
  \left\{
    \begin{aligned}
      dx(t) &=  \left(Ax(t) + I\right)dt \; + \;  \Sigma dB(t)\;,\\
      x(0)&=x_0\;,
    \end{aligned}
  \right.
  \label{linearsde2}
\end{equation}
where $\Sigma \in \R^{m \times p}$ and $B(t)$ is standard
$p$-dimensional Brownian motion.  The following theorem is proved in
the same manner as Theorem \ref{convergence}.

\begin{theorem}  
  The process $\xl(t)=\xl(t,B)$ defined by
  \begin{equation}
    \label{eq:xStationary2}
    \xl(t,B) = \int_{-\infty}^t e^{A(t-s)}I \ ds + \int_{-\infty}^t e^{A(t-s)} \Sigma dB(s)
  \end{equation}
  is a stationary solution to \eqref{linearsde2}. Furthermore given
  any $x_0$, if $x(t,x_0,B)$ is a solution to equation
  \eqref{linearsde2} then $x(t,x_0,B)$ converges to $\xl(t,B)$ as $t
  \rightarrow \infty$ in that
  \begin{equation*}
    \E |x(t,x_0,B) - \xl(t,B)|^2 \rightarrow 0 \quad\text{ as}\quad t \rightarrow \infty\;. 
  \end{equation*}
  \label{convergence2}
\end{theorem}
\begin{proof}
  The proof is identical to that of Theorem \ref{convergence}, except
  that the It\^o Isometry is used to control the expected value of the
  square of the It\^o integral term.
\end{proof}

Since $\xl(t)$ is stationary, the distribution of $\xl(t)$ is
independent of $t$ and invariant under the dynamics of
\eqref{linearsde} (or \eqref{linearsde2}). More precisely, defining the
measure $\mu(A) = \P(\xl(0) \in A)$ for all measurable $A \subset
\R^m$, we see that
\begin{equation*}
  \mu(A)= \int \P(x(t,y,\xi) \in A ) \mu(dy) \;. 
\end{equation*}
Furthermore, $\mu$ characterizes the longtime behavior of the solution
in that the distribution of $x(t,x_0,\xi)$ converges to $\mu$ as $t
\rightarrow \infty$. This follows from $\E|x(t,x_0,\xi)-\xl(t,\xi)|^2
\rightarrow 0$ and the fact that $\mu(A) = \P(\xl(t) \in A)$ for all
$t$.

Thus $\mu$ contains information about the average, long-term behavior
of fluxes and concentrations.  It will be $\mu$, therefore, which we
shall probe in order to gain an understanding of how different
graphical structures and asymptotic limits of biochemical reaction
systems increase, decrease, and otherwise modify the exogenous
fluctuations of biochemical reaction systems.  Throughout the rest of
this paper, it is understood that each mean or variance is computed
with respect to this stationary measure.

\subsection{A General Bound}

We can now prove a simple general bound for the variance of the
concentration of any species in an SSC system in terms of the variance
of the input fluctuations. We assume that the fluctuations, $\xi_t$,
are one-dimensional, stationary, mean zero, and finite variance.  By
taking the expected value in equation \eqref{eq:xStationary} and using
that $\xi_t$ has mean zero one sees that
\begin{equation}
  m_i = I \int_{-\infty}^t e^{A(t-s)}e_1 \cdot e_i  \ ds
\end{equation}
\noindent is the mean of the $i^{th}$ species. 

\begin{theorem} Let $\xl(t)$ be the stationary solution of an SSC
  system with one input, $I$, to a single species, $X_1$, that is
  perturbed by a stationary stochastic process, $\xi_t$, with finite
  variance and mean zero.  Then for each $i$,
\begin{equation*}    
  Var(\xl_i) < \left( \frac{m_i}{I} \right)^2 Var(\xi).
\end{equation*}
\label{declinear}
\end{theorem}
\begin{proof}
  Using Lemma 2.3(c) and the Cauchy-Schwarz inequality gives
  \begin{align*}
    Var(\xl_i(t)) &= \E \left( \int_{-\infty}^t \xi_s e^{A(t-s)}
      e_1 \cdot e_i  \ ds \right)^2\\
    & = \E \left( \int_{-\infty}^t \xi_s \left(e^{A(t-s)} e_1 \cdot
        e_i \right)^{1/2} \left( e^{A(t-s)} e_1 \cdot e_i
      \right)^{1/2} \ ds
    \right)^2\\
    & < \E \left( \int_{-\infty}^t \xi_s^2 e^{A(t-s)} e_1 \cdot e_i \
      ds \right) \left( \int_{-\infty}^t e^{A(t-s)} e_1 \cdot e_i \ ds
    \right)\\
    &= Var(\xi) \left( \int_{-\infty}^t e^{A(t-s)} e_1 \cdot e_i \ ds
    \right)^2 \\
    &= \left( \frac{m_i}{I} \right)^2 Var(\xi).
  \end{align*}
  The strictness of the inequality follows because $\xi_t$ is not a
  constant.
\end{proof}

\vspace{.125in} This simple result is all that we need in this paper.
An analogous proof works in the more general case where there are
inputs to more than one species and any number of the inputs undergo
independent fluctuations.

\section{Chains}

In this section we consider non-reversible chains with mass action
kinetics:
\begin{align}
\begin{array}{ccccccccccccc}
  &I & & k_1 & & k_2 & & k_{m-2} & & k_{m-1}& & k_m &\\ 
0 & \longrightarrow & X_1 & \longrightarrow & X_2 & \longrightarrow & \cdots 
& \longrightarrow & X_{m-1} & \longrightarrow & X_m & \longrightarrow & 0.\\
\end{array}
\label{linearchain}
\end{align}

\noindent Theorem \ref{declinear} allows us to see that variances of
the fluxes of the stationary solution decrease as one proceeds down
the chain.

\begin{theorem} Let the input, $I$, of a non-reversible chain with
  mass action kinetics be perturbed by a stationary stochastic
  process, $\xi_t$, with finite variance and mean zero.  Let $\xl(t)$
  denote the stationary solution for the chain.  Then, for all $i$,
  $Var(k_i \xl_i) < Var(\xi)$ and
  \begin{align}
    Var(k_{i+1} \xl_{i+1}) & \; < \;  Var(k_i \xl_i).
    \label{onelink} 
  \end{align}
\label{ThDecrease}
\end{theorem}
\begin{proof}
  From the remark following Theorem 2.4, we know that the mean, $m_i$, of
  $\xl_i(t)$ is the equilibrium value of $x_i$ for the unperturbed
  problem. For the chain this implies that $m_i = \frac{I}{k_i}$, so
  the bound $Var(k_i \xl_i) < Var(\xi)$ follows immediately from
  Theorem \ref{declinear}.  To prove \eqref{onelink} note that the
  input to $X_2$ is
  \begin{equation*}
    k_1 \xl_1(t) \; = \; I \; + \; (k_1 \xl_1(t) - I) 
  \end{equation*}
  and $k_1 \xl_1(t) - I$ is a stationary stochastic process of mean
  zero and finite variance. Thus, by Theorem \ref{declinear},
  \begin{equation*}
    Var(k_2 \xl_2) \; < \; Var(k_1 \xl_1 - I) \; = \; Var(k_1 \xl_1)\;.
  \end{equation*}
  The input to $X_3$ is $k_2\xl_2(t)$, so repeating this 
  argument down the chain proves \eqref{onelink}.  
\end{proof}

\vspace{.125in}

Note that the variances of the fluxes are strictly decreasing as one
moves down the chain even though the means of the fluxes remain
unchanged (i.e., equal to $I$). The next natural question is how much
do the variances decrease down the chain?  This cannot be answered
without more detailed information about $\xi_t.$ To investigate it, we
will perturb the input $I$ by white noise, $\sigma dB(t),$ which will allow us to use
the It\^o calculus.

\begin{theorem}
  Let $\xl(t)$ be the stationary solution of the linear chain
  \eqref{linearchain} where the input is perturbed by white noise. We
  assume that the rate constants, $k_i$, are distinct.  Then
  \begin{equation}
    Var(\xl_i) = \sigma^2 \sum_{j=1}^i \sum_{r=1}^i p_{ij}p_{ir}\frac{1}{k_j + k_r},
    \label{varcompappendix}
  \end{equation}
  where
  \begin{align}
    p_{ij} &= \left \{ \begin {array}{cc} \left(\prod_{n=1}^{i-1}
          k_n\right)\Big/ \left(\prod_{n=1, l \ne j}^{i} (k_n-k_j)\right) & i \ge j \\ 
        0 & i < j \end{array}.  \right.
  \end{align}
\end{theorem}
\begin{proof}
  The matrix of rate constants, $A$, is given by
  \begin{align*}
    A &= \left[ \begin {array}{cccc}
        -k_{{1}}&0&\dots&0\\\noalign{\medskip}k_{{1}}&-k_{{2}}&\dots&0\\
        \noalign{\medskip}\vdots&\vdots&\ddots&\vdots\\\noalign{\medskip}0&\dots&k_{{m-1}}&-k_{{m}}\end
      {array} \right].
  \end{align*}
  
  \noindent Let $P = \{p_{ij}\}$.  A straightforward calculation shows
  that the $jth$ column of $P$ is the eigenvector of $A$ corresponding
  to eigenvalue $-k_j$. Thus, $D = P^{-1}AP$ is diagonal.  In
  addition, $P$ takes the vector $(1, 1, \cdots, 1)^T$ to the vector
  $(1, 0, \cdots, 0)^T$.  Using these facts, the formula
  \eqref{eq:xStationary2} for $\xl(t)$, and the It\^o Isometry,
  \begin{align*}
    Var(\xl_i) &= \sigma^2 \E \left(\int_{-\infty}^t e^{A(t-s)}e_1
      \cdot e_i dB_s \right)^2\\
    &= \sigma^2 \int_{-\infty}^t \left(Pe^{D(t-s)}P^{-1}e_1 \cdot
      e_i\right)^2 ds \\
    &= \sigma^2 \int_{-\infty}^t
    \left(Pe^{D(t-s)}\left[\begin{array}{c}
          1\\
          \vdots \\
          1
        \end{array}\right] \cdot e_i\right)^2 ds \\ 
    &= \sigma^2 \int_{-\infty}^t \left(P\left[\begin{array}{c}
          e^{-k_1(t-s)}\\ 
          \vdots \\
          e^{-k_m(t-s)}
        \end{array}\right] \cdot e_i\right)^2 ds \\
    &= \sigma^2 \int_{-\infty}^t \left(\sum_{j=1}^i p_{ij} e^{-k_j(t-s)}\right)^2 ds\\
    &= \sigma^2 \sum_{j=1}^i \sum_{r=1}^i p_{ij} p_{ir} \frac{1}{k_j
      + k_r}. 
  \end{align*}
\end{proof}

\vspace{.125in}

We assumed that the $k_i$'s were distinct so that the explicit
formulas above make sense.  It can be shown that the variances of the
concentrations are continuous functions of the rate constants.  This
fact, together with the bound given by \eqref{onelink} allows us to
conclude that formula \eqref{varcompappendix} has finite limits as
various subsets of the $k_i$'s become identical.

\vspace{.125in}

We can use the explicit formula \eqref{varcompappendix} to answer
several natural questions:

\begin{example}\textbf{(Magnitude of decrease)} Theorem
  \ref{ThDecrease} shows that variances of fluxes are strictly
  decreasing as one moves down a chain. To investigate how much they
  decrease, consider the chain \eqref{linearchain} where $m = 2$ and
  the input is perturbed by white noise.  Using
  \eqref{varcompappendix} we see that $Var(k_1\xl_1) =
  \frac{\sigma^2k_1}{2}$ and $Var(k_2\xl_2) =
  \frac{\sigma^2k_1k_2}{2(k_1 + k_2)}$.  Thus,
  \begin{equation*}
    \frac{Var(k_2\xl_2)}{Var(k_1\xl_1)} = \frac{k_2}{k_1 + k_2}\;. 
  \end{equation*}
  This simple example shows that the ratio of successive variances can
  be any number between zero and one.
\end{example}

\begin{example}\textbf{(Long chains)} Assume that $k_i = k$ for
  some fixed $k > 0$ and all $i$.  Taking the limit of
  \eqref{varcompappendix} is difficult.  Instead, since all the $k_i$
  are equal, an induction proof shows that
\begin{equation*}
  \xl_i(t) = \frac{I}{k} + \sigma
  \frac{k^{i-1}}{(i-1)!}\int_{-\infty}^t (t-s)^{i-1}e^{-k(t-s)}dB(s).    
\end{equation*}
Using the It\^o Isometry, it follows that
\begin{equation*}
  Var_{\infty}(\xl_i) \; = \; \sigma^2
  \frac{2(2i-2)!}{4^i(i-1)!^2}\frac{1}{k}\;,   
\end{equation*}
and using Stirling's formula
\begin{equation*}
  Var(k\xl_i) \; \sim \;  \sigma^2
  \frac{k}{2\sqrt{\pi}}\frac{1}{\sqrt{i}} + O(i^{-3/2})\;.    
\end{equation*}
Thus the variances decrease to zero in a regular fashion if all of the rate constants
are the same.
\end{example}

\begin{example} \textbf{(A small rate constant)} Suppose that one
  rate constant, $k_i$, in a chain is very small. Using the explicit
  formula \eqref{varcompappendix}, one can easily compute that
  \begin{align*}
    Var(k_i x_i^*) & \sim \sigma^2 \frac{1}{2}k_i + O(k_i^2),
    \;\;\;\mbox{as}\; k_i \to 0, \\ 
    Var(k_j x_j^*) & \sim \sigma^2 \frac{1}{2} k_i + O(k_i^2), \;\;\;
    \mbox{as} \; k_i \to 0, \; \mbox{for}\; j > i.
  \end{align*}
  Notice that the small rate constant has the effect of significantly
  decreasing the variances of the $i$th and all subsequent fluxes while
  the means of the fluxes remain unchanged. Therefore a small rate
  constant is not ``rate limiting'' but instead is ``variance
  limiting.''
\end{example}

\begin{example}\textbf{(A large rate constant)} Suppose that one
  rate constant, $k_i$, in a chain is very large. Again, using
  \eqref{varcompappendix}, one can compute that
  \begin{align*}
    Var(k_i x_i^*) &\to Var(k_{i-1} x_{i-1}^*), \;\;\; \mbox{as}\; k_i
    \to \infty.
  \end{align*} 
  Furthermore, for all $j > i$,
  \begin{align*}
    Var(k_j x_j^*) & \to Var(k_j \tilde x_j^*), \;\;\; \mbox{as}\; k_i
    \to \infty,
  \end{align*}
  \noindent where $\tilde x_j$ is from the process arising from the following system:
  \begin{figure}[ht]
    \centering
    \begin{picture}(125,30)
      
      \put(-80,20){\makebox(0,0){$I+\sigma dB(t)$}}
      \put(-110,5){\vector(1,0){65}}
      \put(-2,20){\makebox(0,0){$k_1$}}
      \put(-30,5){\makebox(0,0){$\tilde X_1$}}
      \put(-17,5){\vector(1,0){25}}
      \put(23,5){\makebox(0,0){$\tilde X_{2}$}}
      \put(36,5){\vector(1,0){25}}
      \put(49,20){\makebox(0,0){$k_2$}}
      \put(79,5){\makebox(0,0){$\cdots$}}
      \put(130,20){\makebox(0,0){$k_{i-1}$}}
      \put(100,5){\makebox(0,0){$\tilde X_{i-1}$}}
      \put(118,5){\vector(1,0){25}}
      \put(161,5){\makebox(0,0){$\tilde X_{i+1}$}}
      \put(178,5){\vector(1,0){25}}
      \put(189,20){\makebox(0,0){$k_{i+1}$}}
      \put(220,5){\makebox(0,0){$\cdots$}}
    \end{picture}
    \label{newchain}
  \end{figure}
  
  This shows that in the asymptotic limit where $k_i \to \infty$ one
  can replace the original chain by the chain with the substrate $X_i$
  removed. Here we implicitly use the fact that since the kinetics are
  linear and hence the concentrations are Gaussian the statistics are
  determined by the means and variances.

\end{example}
                    
\vspace{.125in}

\section{Side Reaction Systems and Feedback Loops.}

A {\em side reaction system} on a chain is any SSC system that gets
its input from a species on the chain and has output that flows back
into the same species; see Figure \ref{sidechain} below.
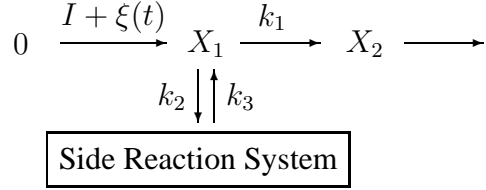
\begin{figure}[ht]
  
  \centering
  \begin{picture}(125,67)
    
    \put(-15,50){\makebox(0,0){$0$}}
    \put(20,60){\makebox(0,0){$I+\xi(t)$}}
    \put(0,50){\vector(1,0){40}}
    \put(52,40){\vector(0,-1){20}}
    \put(58,20){\vector(0,1){20}}
    \put(42,30){\makebox(0,0){$k_2$}}
    \put(68,30){\makebox(0,0){$k_3$}}
    \put(80,60){\makebox(0,0){$k_1$}}
    \put(55,50){\makebox(0,0){$X_1$}}
    \put(37,-10){\makebox(30,30){Side Reaction System}}
    \put(-5,-4){\framebox(114,20){}}
    \put(68,50){\vector(1,0){30}}
    \put(115,50){\makebox(0,0){$X_{2}$}}
    \put(130,50){\vector(1,0){30}}
    
  \end{picture}
  
  \caption{A side reaction on a linear chain}
  \label{sidechain}
\end{figure}

\noindent Note that there must be a species within the side reaction
system whose output flows to $X_1$ with some rate constant, $k_3$.
Define $Y$ to be that species.  The SDE governing the behavior of
$x_1(t)$ is then given by
\begin{equation}
  \frac{d}{dt}x_1(t) = I - k_1x_1(t) - k_2x_1(t) + k_3 y(t)  +  \xi(t).
  \label{whitesidereaction}
\end{equation} 
 If $\tilde x_1$ is the solution to the above system when
there is no side reaction system (i.e. $k_2 = k_3 = 0$), then
\begin{equation}
  \frac{d}{dt}\tilde x_1(t) = I - k_1\tilde x_1(t) + \xi(t).
  \label{nowhitesidereaction}
\end{equation}

\begin{theorem}[Side reactions lower variance] Let $x_1^*$ and $\tilde
  x_1^*$ be the first components of the stationary solutions to
  \eqref{whitesidereaction} and \eqref{nowhitesidereaction},
  respectively, where $\xi(t)$ is a finite variance, mean zero, random
  process or white noise.  Then,
  \begin{equation*}
    Var(k_1 x_1^*) < Var(k_1\tilde x_1^*).
  \end{equation*}
  \label{whitesidethm}
\end{theorem}

\begin{proof} We give the proof in the case where $\xi(t)=\sigma
  dB(t)$ is white noise; the proof in the general case is similar but
  more complicated \cite{anderson1}. Let $z(t)= \E (k_1x_1(t) - I)^2$
  and $\tilde z(t)= \E (k_1 \tilde x_1(t) - I)^2,$ where $x_1(t)$ and
  $\tilde x_1(t)$ are solutions of \eqref{whitesidereaction} and
  \eqref{nowhitesidereaction}.  By theorem \ref{convergence2} $z(t)$
  and $\tilde z(t)$ converge to $Var(k_1 \xl_1)$ and $Var(k_1 \tilde
  x_1^*)$, respectively.  We will prove the theorem by comparing the
  differential equations satisfied by $z(t)$ and $\tilde z(t)$.
  
  By using Kolmogorov's backward equation \cite{oxy}, we see that $\tilde z(t)$ satisfies
  \begin{equation}
    \tilde z'(t) = -2k_1 \tilde z(t) + k_1^2 \sigma^2.
    \label{tildeeq}
  \end{equation}
  Therefore, $Var(k_1 \tilde x_1^*) = k_1 \sigma^2/2$ because
  $Var(k_1\tilde x_1^*)$ is the equilibrium value of \eqref{tildeeq}.
  Similarly, $z(t)$ satisfies 
  \begin{equation*}
    z'(t) = -2k_1z(t) + k_1^2 \sigma^2 + 2k_1 \mathbb E(k_1x_1(t) - I)(k_3 y(t) - k_2 x_1(t)),
  \end{equation*}
  and so, by Theorem \ref{convergence2}, 
  \begin{equation*}
    Var(k_1 \xl_1) = \frac{k_1 \sigma^2}{2} + \frac{1}{2} \E(k_1 \xl_1 - I)(k_3 y^* - k_2 \xl_1).
  \end{equation*}
  Thus, to complete the proof we need only show that $\mathbb E (k_1
  \xl_1 - I)(k_3 y^* - k_2 \xl_1) < 0$.  The remark following Theorem
  2.4 implies that $\mathbb E \xl_1 = I/k_1$ and $\mathbb E k_3 y^* =
  \mathbb E k_2 \xl_1$.  Therefore, $\mathbb E k_3 y^* = \frac{k_2
    I}{k_1},$ and
  \begin{equation*}
    \E  (k_1 \xl_1 - I )(k_3 y^* - k_2\xl_1) = \frac{k_1}{k_2} \E \left( k_2
      \xl_1 k_3 y^* - k_2^2 {\xl_1}^2 \right). 
  \end{equation*}
  By Theorem \ref{declinear} \;  $\E \left( k_3 y^*\right) ^2 < \mathbb E
  \left( k_2 \xl_1 \right)^2$, so
  
  \begin{align*}
    \left| \E \left( k_2 \xl_1 k_3 y^*\right) \right| & \le \left(\E \
      k_2^2 {\xl_1}^2 \right)^{\frac{1}{2}} \left( \E \ k_3^2
      {y^*}^2\right)^{\frac{1}{2}} \\
    & < \E \ k_2^2 {\xl_1}^2.
  \end{align*}
  Thus, $\E(k_1 \xl_1 - I)(k_3 y^* - k_2\xl_1) < 0$, as desired.
\end{proof}

\lskip

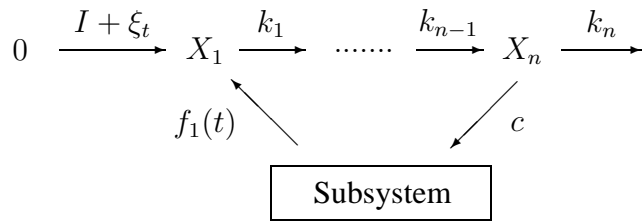
\begin{figure}[ht]
  \centering
  \begin{picture}(80,79)
    \put(-95,60){\makebox(0,0){$0$}}
    \put(-60,70){\makebox(0,0){$I+\xi_t$}}
    \put(-80,60){\vector(1,0){40}}
    \put(0,70){\makebox(0,0){$k_1$}}
    \put(-25,60){\makebox(0,0){$X_1$}}
    \put(27,-10){\makebox(30,30){Subsystem}}
    \put(0,-4){\framebox(83,20){}}
    \put(-12,60){\vector(1,0){25}}
    \put(95,60){\makebox(0,0){$X_{n}$}}
    \put(110,60){\vector(1,0){30}}
    \put(35,60){\makebox(0,0){$.......$}}
    \put(55,60){\vector(1,0){25}}
    \put(67,70){\makebox(0,0){$k_{n-1}$}}
    \put(125,70){\makebox(0,0){$k_{n}$}}
    \put(93,48){\vector(-1,-1){25}}
    \put(10,24){\vector(-1,1){25}}
    \put(93,32){\makebox(0,0){$c$}}
    \put(-25,32){\makebox(0,0){$f_1(t)$}}
  \end{picture}
  \caption{A chain with a feedback loop}
  \label{feedbackloop}
\end{figure}

A {\em feedback loop} on a chain is an SSC system together with an
input from one species on the chain, $X_n$, and an output to an
earlier species, $X_1$; see Figure \ref{feedbackloop}.

\begin{theorem}
  Let $\tilde x(t)$ be the vector of species concentrations for the
  chain \eqref{linearchain} and let $x(t)$ be the vector of species
  concentrations for the chain with feedback loop (Figure
  \ref{feedbackloop}), where $\xi(t)$ is a finite variance, mean zero,
  random process or white noise. Then,
  \begin{equation*}
    Var(k_n \xl_n) < Var(k_n {\tilde x}_n^*).
  \end{equation*}
  \label{sideloop}
\end{theorem}

\begin{proof} Let $\{V_i\}$ be the substrates and $B$ be the matrix of
  rate constants of the SSC subsystem in Figure \ref{feedbackloop}.
  We suppose that $V_j$ is the species which gives input to $X_1$ with
  rate constant $\alpha$.  Then the input to $X_1$ from the feedback
  loop is
  \begin{equation*}
    f_1(t)= \alpha e^{Bt}v(0) \cdot e_j +  \alpha c \int_0^t
    e^{B(t-s)}x_n(s) \cdot e_j \ ds, 
  \end{equation*}
  which depends explicitly only on $x_n$.  If we let $R(t) =
  k_{n-1}x_{n-1}(t)$ then the differential equation for $x_n(t)$ is
  $\dot x_n(t) = R(t) - cx_n(t) - k_nx_n(t)$.
  
  \begin{figure}[ht]
    
    \begin{picture}(80,130)(-160,5)
      
      \put(-105,120){\makebox(0,0){$0$}}
      \put(-70,130){\makebox(0,0){$I+\xi_t$}}
      \put(-90,120){\vector(1,0){40}}
      \put(-10,130){\makebox(0,0){$k_1$}}
      \put(-35,120){\makebox(0,0){$X_1$}}
      \put(-22,120){\vector(1,0){25}}
      \put(85,120){\makebox(0,0){$X_{n}$}}
      \put(100,120){\vector(1,0){30}}
      \put(25,120){\makebox(0,0){$.......$}}
      \put(45,120){\vector(1,0){25}}
      \put(57,130){\makebox(0,0){$k_{n-1}$}}
      \put(115,130){\makebox(0,0){$k_{n}$}}
      \put(75,108){\vector(-1,-1){65}}
      \put(35,80){\makebox(0,0){$c$}}
      \put(-72,11){\framebox(83,20){}}
      \put(-47,5){\makebox(30,30){Subsystem}}
      \put(20,19){\vector(1,0){40}}
      \put(40,5){\makebox(0,0){$f_1(t)$}}  
      \put(71,19){\makebox(0,0){$Y_1$}}
      \put(79,19){\vector(1,0){20}}
      \put(89,5){\makebox(0,0){$k_1$}}
      \put(109,19){\makebox(0,0){$Y_2$}}
      \put(119,19){\vector(1,0){20}}
      \put(129,5){\makebox(0,0){$k_2$}}
      \put(155,19){\makebox(0,0){$\cdots$}}
      \put(182,19){\makebox(0,0){$Y_{n-1}$}}
      \put(174,30){\vector(-1,1){77}}
      \put(140,80){\makebox(0,0){$k_{n-1}$}}
      
    \end{picture}
    
    \caption{A chain with a side reaction system}
    \label{feedbackasside}
  \end{figure}
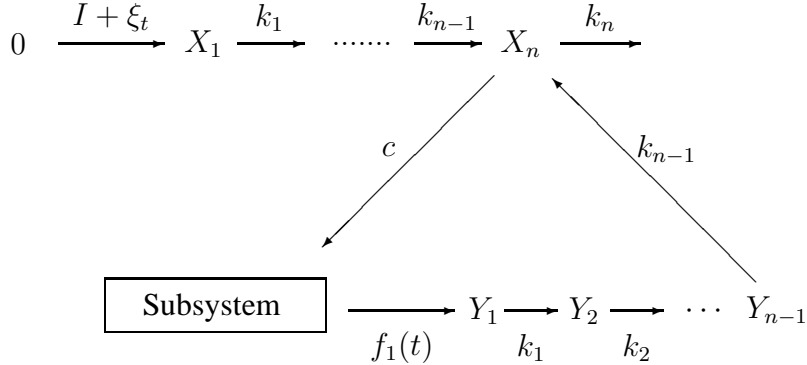
  Consider the chain with side reaction system given in Figure
  \ref{feedbackasside} where the subsystem is the same as in Figure
  \ref{feedbackloop} and the flux to $Y_1$ comes from $V_j$ with rate
  constant $\alpha$.  Let $Q(t) = k_{n-1}x_{n-1}(t)$ and $P(t) =
  k_{n-1} y_{n-1}(t)$ be the inputs to $X_n$ in Figure
  \ref{feedbackasside}.  Since the input to the Y-chain is $f_1(t)$
  and the rate constants for the two chains are the same, $R(t) = Q(t)
  + P(t)$ because the differential equations are linear.  Thus, the
  differential equation governing $x_n(t)$ in Figure
  \ref{feedbackloop} is the same as the differential equation
  governing $x_n(t)$ in Figure \ref{feedbackasside}.  Since the system
  in Figure \ref{feedbackasside} is a chain with a side reaction
  system, the result follows from Theorem \ref{whitesidethm}.
\end{proof}

\section{One large rate constant in a general SSC system}

We now consider a general weakly reversible SSC system with input and
characterize the effect of a large rate constant.

\begin{theorem}
  Suppose that independent white noise processes perturb the inputs to
  a weakly reversible SSC system with $m$ substrates.  Let $X_a$ be a
  particular substrate and suppose that the rate constant $L$ for one
  flux out of $X_a$ to another complex (possibly the zero complex) is
  large.  Then,
  \begin{equation}  
    Var(x_a^*) \sim O\left(\frac{1}{L}\right)  \;\;\; \mbox{as} \; L \to \infty.
  \end{equation}
  \label{asymptotics}
\end{theorem}
\begin{proof}
  We will assume that one of the perturbed inputs goes directly to
  $X_a$.  The proof of the general case is similar.  The stochastic
  differential equation governing $x_a(t)$ is given by
  \begin{equation}
    dx_a(t) = \left(C + \sum_{i=1}^m c_i x_i(t) - (L + K)x_a(t)\right)dt
    + \sigma dB(t)\;, 
    \label{Lsde}
  \end{equation} 
  where $L+K >0$ is equal to the sum of all the rate constants for
  reactions leaving $X_a$, $C > 0$ is the input flow to $X_a$ from the
  zero complex, $\sigma > 0$, and $c_i \ge 0$ is the rate constant
  associated with the reaction $X_i \rightarrow X_a$.  Solving
  \eqref{Lsde} for $x_a^*$ in terms of the $x_i^*$ and using the It\^o
  Isometry, one can easily bound $Var(x_a^*)$,
  \begin{equation*}
    Var(x_a^*)\; \le \; \frac{\beta}{2(L+K)} +  \beta
    \sum_{i=1}^m \frac{c_i^2Var(x_i^*)}{(L+K)^2}, 
  \end{equation*}
  for some constant $\beta.$ To complete the proof we will show that
  $Var(x_i^*) \leq O(L)$.
  
  \vspace{.125in}
  
  Let $A$ be the matrix of rate constants for the SSC system. Using
  the formula \eqref{eq:xStationary2} for the stationary solution and
  the Ito Isometry, one easily calculates:
  \begin{equation}
    Var(x_i^*) = \sigma^2 \int_{-\infty}^t (e^{A(t-s)} e \cdot e_i)^2  \ ds,
    \label{matrixbound}  
  \end{equation}
  for some vector $e$.  By Lemma \ref{linearstable}(b) we know that
  the real parts of the eigenvalues of $A$, $\{\lambda_i\}$, are
  strictly negative; let $\lambda = \inf{\{|\lambda_i|\}}$.  There
  exist positive constants $c$ and $M$ so that for all $t-s>0$, we
  have $\|e^{A(t-s)}\| \le c e^{-M\lambda(t-s)}$.  Using this
  inequality in \eqref{matrixbound}, we have
  
  \begin{equation*}
    Var(x_i^*) \le \frac{\sigma^2 c^2|e|}{2 M}\frac{1}{\lambda}.
  \end{equation*}
  In Appendix A we prove that $\lambda \ge O(1/L)$, so $Var(x_i^*)
  \leq O(L)$, which concludes the proof.
\end{proof}

\noindent \textbf{ Example 6.1 (A side chain with a large rate constant)}
To illustrate the theorem, we consider the linear chain with a side
reaction given in the diagram in the Introduction. As the rate
constant $L$ becomes large, Theorem \ref{asymptotics} tells us that
$Var(x_2^*) \leq O(1/L)$. Therefore the flux out of $X_2$ down the
chain has variance $Var(k_2x_2^*) \leq O(1/L)$.  By Theorem
\ref{declinear},
\begin{equation*}
  Var(k_ix_i^*) \leq Var(k_2x_2^*) \le O(1/L) \hspace{.2in} \mbox{for
  all $i \geq 2.$} 
\end{equation*}
Thus, for all $i \geq 2$, the means of the fluxes remain equal to $I$,
while the variances of the fluxes go to zero as $L \to \infty$.

\section{Application to Methionine Metabolism.}

The actual biochemical systems involved in cell metabolism are much
more complicated and more difficult to analyze than the single species
systems considered in the previous sections. Consider, for example the
diagram in Figure 6.1 that shows the methionine cycle and part of the
folate cycle. Firstly, most reactions have two or more substrates and
many enzymes are inhibited by the products of the reactions they
catalyze. Thus, the kinetics will be highly nonlinear. Secondly, many
reactions are catalyzed by two different enzymes that have very
different properties.  Thirdly, some substrates inhibit or activate
distant enzymes in the reaction diagram (red arrows in the diagram).
These long-range interactions make it virtually impossible to intuit
the emergent properties of the network by tracing influences from
point to point.
\begin{center}
  \includegraphics[width= 4.5in]{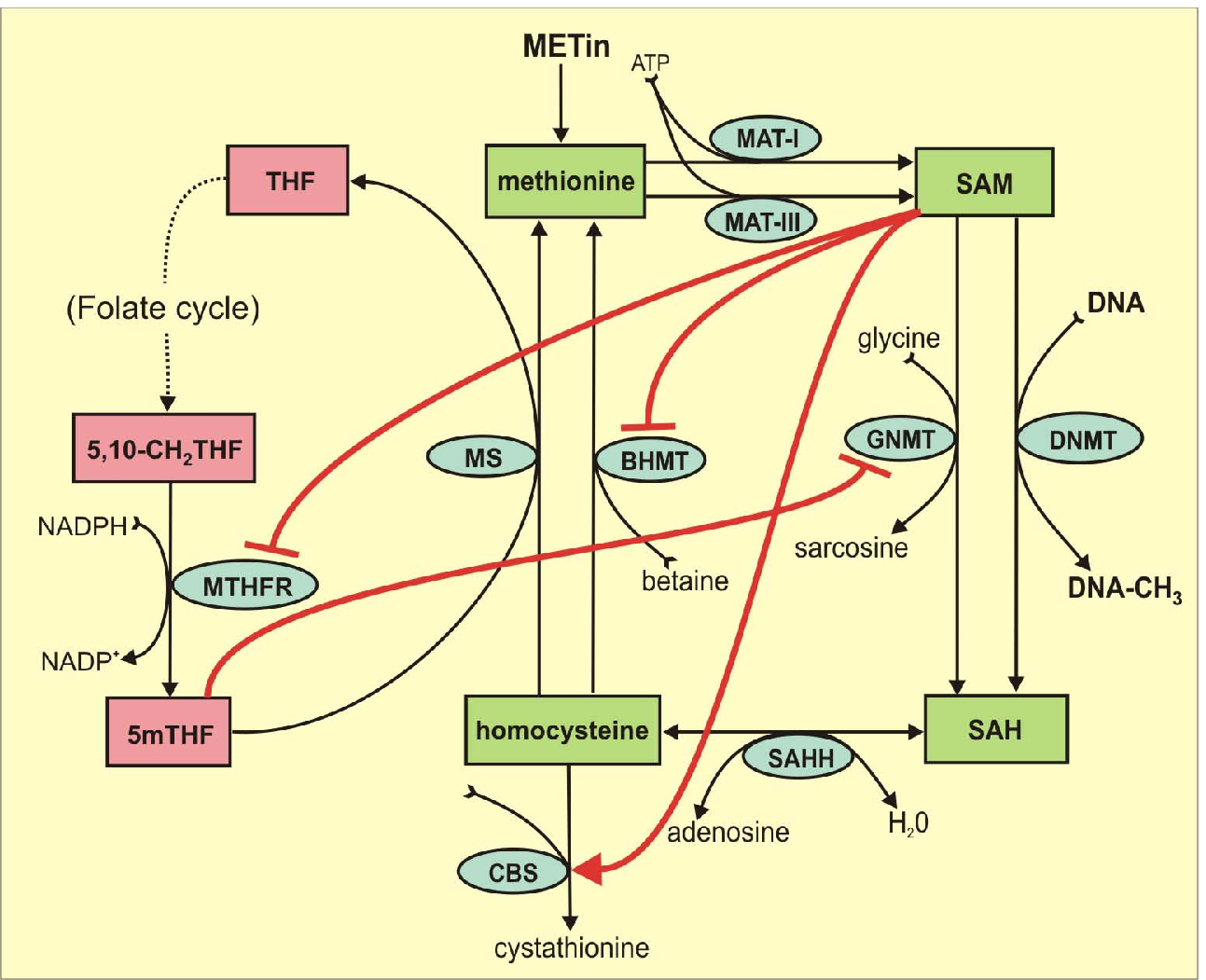}
\end{center}

{\small { \bf Figure 6.1. Methionine Metabolism.} Substrates of the
  methionine cycle and (part of) the folate cycle are shown in green
  and red rectangles, respectively. Enzyme acronyms are in ellipses.
  Long-range interactions are shown by red curves with the arrow
  indicating activation and the bars indicating inhibition. SAM,
  s-adenosyl-methionine, activates CBS and inhibits BHMT and MTHFR,
  while 5mTHF, 5-methyltetrahydrofolate, inhibits GNMT. }

\vspace{.125in} Epidemiological evidence correlates changes in folate
and methionine metabolism to serious human health consequences
(cancer, heart disease, depression) and there are several important
public health issues involved in folate supplementation as currently
practiced in the United States and Canada. Thus, this part of cell
metabolism has been the object of numerous experimental studies and
several modeling studies
\cite{martinov1}\cite{reed}\cite{nijhout1}\cite{martinov2}\cite{nijhout2}.
Our purpose here is simply to illustrate how fluctuation analysis can
be used to understand such a complex system.

The velocities of the individual reactions in the methionine cycle
\cite{nijhout2} are typically highly nonlinear functions that depend
on the concentrations of several substrates. For example, the velocity
of the GNMT reaction, $V_{GNMT}$, depends on $SAM$, on $SAH$ because
of product inhibition, and on $5mTHF$ because of a long-range
interaction. Because of the complexity and the nonlinearities, a
rigorous mathematical analysis of this system is beyond current
mathematical techniques. Even proving the existence and uniqueness of
a stationary measure is a delicate issue.  Nevertheless, we have
investigated this question by numerical computation in the case where
the methionine input ($METin$) is an Ornstein-Uhlenbeck process with
mean 100 $\mu$M/hr and standard deviation 30 (variance = 900). We
found that the joint distribution of the substrates does indeed
stabilize as time gets large, and thus for each concentration or flux,
X, we can compute the ratio
\[ r \; = \; \frac{ \mbox{Variance}(X)}{\mbox{Variance}(METin)}, \]
which tells us how much $X$ varies compared to the variance of the
input.  Table 1 shows the values of $r$ for two substrates and two
fluxes in the case where all the long-range interactions are present
(regulated) and the case where the long-range interaction are absent
(unregulated). Methionine is quite stable in both cases but is more
stable in the regulated case. $SAM$ is much less stable than
methionine, which agrees with what is seen experimentally. Notice that
in the unregulated case, the variances of $V_{GNMT}$ and $V_{DNMT}$
are similar, but in the regulated case the variance of $V_{GNMT}$
doubles and the variance of $V_{DNMT}$ becomes exceptionally small.
There are good biological reasons why one would want the DNA
methylation rate to be stabilized against fluctuations in methionine
input. Thus, fluctuation analysis shows that this stabilization is
achieved by the long-range interactions. We have also computed the
values or $r$ in all the intermediate cases where some but not all of
the long-range interaction are present and this has enabled us to
quantify each of their effects and propose an evolutionary scenario
\cite{nijhout2}.

 \begin{center}
\begin{tabular}{| c | c | }
\multicolumn{2}{l}{\textbf{Table 1. Values of r}}\\
\hline Methionine & r   \\
\hline  regulated & .064 \\
\hline unregulated & .082 \\ 
\hline 
\hline  SAM & r   \\
\hline  regulated & .22 \\
\hline unregulated & 1.23 \\ 
\hline 
\hline  $V_{GNMT}$ & r   \\
\hline  regulated & .15 \\
\hline unregulated & 0.079 \\ 
\hline 
\hline  $V_{DNMT}$ & r   \\
\hline  regulated & 0.007 \\
\hline unregulated & 0.09 \\ 
\hline
\end{tabular}
\end{center}

In liver cells the reaction from methionine to $SAM$ is catalyzed by
two isoforms of the same enzyme, $MAT\hspace{-.05in}-\hspace{-.05in}I$
and $MAT\hspace{-.05in}-\hspace{-.05in}III$, that have very different
properties \cite{mato}.  $MAT\hspace{-.05in}-\hspace{-.05in}I$ is
inhibited by $SAM$ and $MAT\hspace{-.05in}-\hspace{-.05in}III$ is
activated by $SAM$, and it has been proposed that it is this unusual
combination that stabilizes the methionine concentration.  To test
this, we recomputed $r$ after eliminating the
$MAT\hspace{-.05in}-\hspace{-.05in}III$ reaction and raising the
$V_{max}$ of the $MAT\hspace{-.05in}-\hspace{-.05in}I$ reaction so
that it carried the same flux previously carried by both. The values
in Table 2 show conclusively that, indeed, the presence of the
$MAT\hspace{-.05in}-\hspace{-.05in}III$ reaction somewhat destabilizes
$SAM$ but greatly increases the stability of the methionine
concentration.

 \begin{center}
\begin{tabular}{| c | c | }
\multicolumn{2}{l}{\textbf{Table 2. Values of r}}\\
\hline Methionine & r   \\
\hline with $MAT\hspace{-.05in}-\hspace{-.05in}III$ & 0.06 \\
\hline no $MAT\hspace{-.05in}-\hspace{-.05in}III$ & 0.16 \\ 
\hline 
\hline  SAM & r   \\
\hline   with $MAT\hspace{-.05in}-\hspace{-.05in}III$& 0.22 \\
\hline  no $MAT\hspace{-.05in}-\hspace{-.05in}III$& 0.17 \\
\hline 
\end{tabular}
\end{center}

\vspace{.125in}
\section{Discussion.} 

In Sections 2-5, we developed the theory of propagation of
fluctuations for the special case of linear SSC networks and proved
theorems relating variances to network structure.  Variances decrease
down a chain and the presence of side reactions and feedback loops
always lowers the variances further down the chain. These results are
very general in that they hold independent of the choice of rate
constants. It is tempting to speculate that biochemical systems
evolved to be as complicated as they are partly because of the
homeostasis of exit fluxes achieved by having many intermediate steps.
We also showed how the large size of a single rate constant affects
variances. It is known that most of these results generalize to
non-linear SSC networks with restrictions on the nature of the
nonlinearity \cite{anderson2}. It remains to be seen whether they
generalize to networks in which complexes contain more than one
species. In these highly non-linear contexts, a fundamental
mathematical issue is the proof of the existence of a stationary
measure.

A reasonable concern with the idealized models in Sections 2-5 is
that, under the influence of the fluctuations, the concentrations can
become negative. By modifying the forcing processes appropriately this
could have been avoided. However, this would complicate the analysis
and prevent us from obtaining explicit formulae and straightforward
bounds. Since our goal with these idealized models is to build
intuition and develop general principles, we have purposely avoided
complicating the analysis.

In Section 6 we showed how the ideas of fluctuation theory could be
used to investigate a network of biological interest, the methionine
cycle. It is reasonable to ask whether methionine input actually
fluctuates randomly and if so what are the properties of the
fluctuations. There are really two answers. The input to the
methionine pool in liver cells is certainly continually varying. There
are large deviations on the time scale of hours depending on the times
and content of meals. Methionine is always being used for protein
synthesis and is being made available by protein catabolism, two
processes that are themselves variable and not always in balance. The
methionine available for input to the methionine cycle is also
affected by the use of methionine in other metabolic reactions.
Finally, all these processes are affected by the time-varying
regulation of the genes that code for the various enzymes. Thus, the
first answer is that we don't know how methionine input varies but it
certainly fluctuates with standard deviations of the order of 30-50
$\mu$M/hr on the time scale of hours and with smaller standard
deviations on the time scales of minutes and seconds. The second
answer is that it doesn't matter.  We are using the fluctuations in
methionine input as a probe of the dynamical properties of the system
away from equilibrium. Of course, we need to be sure that the
properties we find do not depend on the detailed properties of the
noise.

For simplicity of exposition, we have discussed the special case where
a single input to a biochemical system is varied.  The same ideas can
be used to introduce fluctuations in a concentration, a flux, or in
several places, and then study how the fluctuations propagate
throughout the system. Understanding the consequences of fluctuations
in kinetic parameters is also important because kinetic parameters
depend on enzyme concentrations and other properties that are variable
and themselves dependent on time-varying genetic regulation. Analyzing
this case requires some technical extensions of this work.

In the Introduction we referred to ``intrinsic stochasticity'' in
contrast to the external stochastic forcing that we consider. It would
be interesting to consider models with both forms of stochasticity,
and indeed both surely arise in gene networks. In gene networks that
are coupled to biochemical networks, the intrinsic stochasticity at
the gene and gene regulation level can be viewed as external
stochastic forcing to the biochemical level. Therefore, both types of
questions and analyses will be necessary to gain full understanding of
real biological networks.

\vspace{.25in}

\noindent \textbf{\Large{Appendix A}}

\lskip

\noindent We derive the bound used in Theorem \ref{asymptotics}.
There are two cases which need consideration:
\begin{enumerate}
\item  The flux out of $X_a$ with rate constant $L$ goes to another
  species.  This case is handled in Theorem A.1 below.  
\item  The flux out of $X_a$ with rate constant $L$ leaves the
  system.  The proof of the result in this case is similar to the proof of
  the theorem below and so the details are omitted.
\end{enumerate}

\lskip

\noindent \textbf{Theorem A.1} {\it Let $A = \{a_{ij}\}$ be an $n
  \times n$ matrix with the following properties:
\begin{enumerate}[(1)]
\item For each $i$, $a_{ii} < 0$ and $|a_{ii}| \ge \sum_{j \ne i}^n
  |a_{ji}|$.
\item $a_{11} = -L + \alpha_{11}$ and $a_{21} = L + \alpha_{21}$ for
  some $\alpha_{11} < 0$ and $\alpha_{21} \in \mathbb R$.
\item For every $L > 0$, the real parts of the eigenvalues of $A$ are
  all strictly negative.
\end{enumerate}
Denote the eigenvalues of $A$ by $\{\lambda_i\}$ and
let $\lambda = \inf{\{|Re(\lambda_i)|\}}$.  Then 
\begin{equation*}
\lambda \ge O(1/L), \hbox{ as } L \to \infty.
\end{equation*}}
\begin{proof} Let $B = \frac{1}{L}A$.  The eigenvalues of $B$ are
  $\{\frac{1}{L}e_i : e_i \hbox{ is an eigenvalue of } A\}$.  We will
  use the characteristic polynomials of $A$ and $B$ to show that the
  magnitude of the real parts of the eigenvalues of $B$ are no smaller
  than $O(1/L^2)$, which implies our result.
  
  Because $L$ only appears in the first column of $A$, all $O(1)$
  terms of $B$ occur in the first column.  Expanding the determinant
  of $B$ by cofactor expansion along the first column then shows that
  $det(B)$ must be of order $O(1/L^n)$ or $O(1/L^{n-1})$.  Similarly,
  the cofactors of $B$ must be of order $O(1/L^{n-1})$ or
  $O(1/L^{n-2})$.  Therefore, computing the inverse of $B$ (which
  exists by assumption (3) above) by cofactors, we see that the
  possible order of the entries of $B^{-1}$ are $1$, $L$, and $L^2$.
  Therefore, $\|B^{-1}\| \leq O(L^2)$.
  
  One may view $B$ as a $1/L$ matrix perturbation of the matrix $C =
  \{c_{ij}\}$, where $c_{11} = -1$, $c_{21} = 1$, and $c_{ij} = 0$ for
  all other entries.  Therefore, each eigenvalue, $\rho$, of $B$ is an
  analytic functions of $1/L$:
  \begin{equation}
    \rho = \rho_0 + \frac{1}{L}\rho_1 + \frac{1}{L^2}\rho_2 +
    O\left(\frac{1}{L^3}\right), 
    \label{eigen_expansion}  
  \end{equation}
  where $\rho_0$ is $-1$ or $0$.  If $\rho_0 =-1$ there is nothing to
  prove; so we assume $\rho_0 = 0.$ If $\rho_1 = \rho_2 = 0$ then
  $\rho = O(1/L^3)$.  However, this would imply that $O(1/\rho) =
  O(L^3)$.  Since $1/\rho$ is an eigenvalue of $B^{-1}$, this would
  contradict the norm bound for $B^{-1}$, above.  Thus $\rho_1$ and
  $\rho_2$ can not both be zero.  It remains to be shown that the
  leading order term in equation \eqref{eigen_expansion} can not be
  purely imaginary.  We will do this through asymptotic matching.  
  
  Consider two different formulations for the characteristic
  polynomial of $A$, $p_A(x)$:
  \begin{align}
    p_A(x) &= det(xI_n - A)\\
    & = x^n + Lu(x) + v(x) \label{char_poly1}\\
    \begin{split}
      & = x^n + c_{1,n-1}L x^{n-1} + c_{0,n-1}x^{n-1} + \cdots +
      c_{1,2}Lx^2 + c_{0,2}x^2 \\
      & \;\;\; + c_{1,1}Lx + c_{0,1}x +c_{1,0}L + c_{0,0},
    \end{split}
    \label{char_poly2}
  \end{align}
  where $u(x)$ and $v(x)$ are polynomials of degree $n-1$ that are
  independent of $L$, and $c_{i,j} \in \mathbb R$ for $i = 1,2$ and $j
  = 1,..n-1$ ($i$ gives the power of $L$ and $j$ gives the power of
  $x$ for the term $c_{ij}L^i x^j$).  We note that we can not have
  $c_{1,0} = c_{0,0} = 0$, for then there would be a zero eigenvalue,
  which would contradict assumption (3).  

  To show that the leading order term in equation
  \eqref{eigen_expansion} is not purely imaginary we will consider two
  cases: $\rho_1 = 0$ and $\rho_1 \ne 0$.  We begin by supposing
  $\rho_1 =0$ and $\rho_2 \ne 0$.  Then $\rho = O(1/L^2)$ and there is
  a solution to \eqref{char_poly2} which is $O\left(1/L\right)$.
  Putting $x = \rho_2/L$ into \eqref{char_poly2} and setting the
  equation equal to zero gives us:
  \begin{equation*}
    O\left(\frac{1}{L^3}\right) + \frac{c_{1,2}\rho_2^2}{L} +
    \frac{c_{0,2}\rho_2^2}{L^2} + 
    c_{1,1}\rho_2 + \frac{c_{0,1}\rho_2}{L} +c_{1,0}L + c_{0,0} = 0.
  \end{equation*} 
  Matching like terms in $L$ tells us that $c_{1,0}=0$, $c_{0,0} \ne
  0$, and $c_{1,1} \ne 0$.  Solving for $\rho_2$ gives us $\rho_2 =
  -c_{0,0}/c_{1,1} \in \mathbb R$.  Therefore, $\rho_2$ has a nonzero
  real part.
  
  We now suppose that $\rho_1 \ne 0$.  Because finding an $O(1/L)$
  solution to equation \eqref{eigen_expansion} is equivalent to
  finding an $O(1)$ solution to \eqref{char_poly1}, $\rho_1$ must
  satisfy $u(\rho_1)=0$.  Let $D(x) = xI_n - A$.  Then $u(x) =
  D(x)_{11} + D(x)_{21}$, where $D(x)_{ij}$ is the $i,j^{th}$ cofactor
  of $D(x)$. $D(x)_{11}$ and $D(x)_{21}$ differ only in the first row,
  so we may combine the determinants by adding the first two rows.  We
  conclude that
  \begin{equation*}
    u(x) = \left|\begin{array}{ccccc}
        -a_{22} -a_{12} + x & - a_{13} - a_{23} & -a_{14} - a_{24} &
        \cdots &  -a_{1n}-a_{2n}\\
        -a_{32} & -a_{33} + x & -a_{34} & \cdots &  - a_{3n}\\
        - a_{42} & - a_{43}  & -a_{44} + x & \cdots &  - a_{4n}\\
        \vdots & \vdots & & \ddots & \vdots\\
        - a_{n2} & - a_{n3}  & - a_{n4} & \cdots &   -a_{nn} + x\\
      \end{array}\right|.
  \end{equation*}
  Solving $u(x) = 0$ for non-zero solutions is therefore
  equivalent to finding the non-zero eigenvalues of the matrix
  \begin{equation*}
    \tilde A = \left[\begin{array}{ccccc}
        a_{22} + a_{12} & a_{13} + a_{23} & a_{14} + a_{24} &
        \cdots &  a_{1n}+a_{2n}\\
        a_{32} & a_{33} & a_{34} & \cdots &  a_{3n}\\
        a_{42} & a_{43}  & a_{44} & \cdots &  a_{4n}\\
        \vdots & \vdots & & \ddots & \vdots\\
        a_{n2} & a_{n3}  & a_{n4} & \cdots &   a_{nn} \\
      \end{array}\right].
  \end{equation*}
  By assumption (1), the diagonal entries of $\tilde A$ are
  non-positive and have magnitudes that are greater than or equal to
  the sums of the magnitudes of all the other terms in that column.
  Therefore, Gershgorin's Theorem says that the non-zero eigenvalues
  of $\tilde A$, and hence the non-zero solutions of $u(x)=0$, have
  strictly negative real part.  Thus, $Re(\rho_1) \ne 0$.  This
  completes the proof.
\end{proof}

If the flux out of $X_a$ with rate constant $L$ leaves the system, the
only change in the statement of the above theorem is that $a_{21}$ is
independent of $L$.  The proof is identical except that $u(x) =
D(x)_{11}$ and so we no longer have to add two determinants together
to simplify $u(x)$.

\vspace{.25in}
\begin{center}
{\bf \large Acknowledgements}
\end{center}
This research was supported by NIH grant R01-CA105437, NSF grants
DMS0109872, DMS0449910, and the Alfred P. Sloan Foundation
(Mattingly).


\pagestyle{plain}

\bibliographystyle{amsplain}

\end{document}